\documentclass[twoside,10pt]{article}
\usepackage{amssymb}

\newtheorem{thm}{Theorem}

\def \k {\kern 10 pt}
\def \n {\noindent}

\title{\bf THE GAUSS-BONNET THEOREM FOR VECTOR BUNDLES}
\begin{document}
\date{}
\maketitle

\kern - 50 pt

\begin{center} \author{ Denis Bell\footnote{Research partially 
supported by
NSF grant DMS-9703852}
\\Department of Mathematics, University of North Florida
\\4567 St. Johns Bluff Road South,Jacksonville, FL 32224,
U. S. A.
\\email: dbell@unf.edu}

\end{center}

\kern 20 pt

\centerline {This paper is dedicated to the memory of
Philip Bell.}

\kern 50 pt

\n {\bf Abstract}. We give a short proof of the Gauss-Bonnet theorem for
a  real oriented Riemannian
vector bundle $E$ of even rank over a
closed compact orientable manifold $M$. This theorem reduces to the
classical
Gauss-Bonnet-Chern theorem in the special case when $M$ is a Riemannian
manifold and $E$ is the tangent bundle of $M$ endowed with the 
Levi-Civita connection. The proof is based on an explicit geometric
construction of the Thom class for 2-plane bundles.

\kern 20 pt

\n Mathematics Subject Classification (2000): 58A10, 53C05, 57R20.

\n Key words: Gauss-Bonnet formula, Thom class, Euler class, metric
connection.

\vfil\break

\n {\bf 1. Introduction}

\kern 2 pt

 The classical Gauss-Bonnet theorem gives a remarkable relationship between
the topology and the geometry of a compact orientable surface in {\bf
R}$^3$. In
1944, Chern generalized this theorem to all
even-dimensional compact orientable manifolds, proving what is now known
as the 
Gauss-Bonnet-Chern theorem. The Gauss-Bonnet-Chern theorem occupies a
central place
in differential geometry, opening the way to such important developments
as the
theory of characteristic classes and index theory. It continues to be a
subject of
current research, and new proofs of it have appeared within the last
decade, e.g. Hsu
[4] and Rosenberg [6]. The most recent proofs are based upon the heat
kernel method pioneered by Mckean-Singer and Patodi, an approach that
was reinvigorated by Bismut's 
 stochastic proof [1] of the index theorem for the Dirac operator on
spin bundles. 

In this paper we give a short proof of a more general theorem. In
contrast to the
heat
kernel approach,
the argument presented here avoids hard analytical estimates and 
provides direct
insight into the structure of the integrand in the Gauss-Bonnet
formula. We prove the following:

 \begin{thm}  Let $E$ denote a real
orientable Riemannian vector bundle of even rank $n=2k$ over a
closed compact  manifold $M$ of dimension $n$. Then
$$\Big ({-1\over 2\pi}\Big )^{k}
\int_M Pf(\Omega) = \chi (E) \eqno (1)$$  
where $Pf$ denotes the Pfaffian function acting on square skew-symmetric
matrices,
$\Omega$ the curvature 2-form with respect to any metric
connection on $E$,  and $\chi (E)$ the Euler characteristic of $E$.
\end{thm}

This result reduces to the Gauss-Bonnet-Chern theorem in the case
when
$M$ is a Riemannian manifold and $E$ is the tangent bundle of $M$ endowed
with the Levi-Civita connection.

Section 2 contains background
information. In particular, the
curvature of the vector bundle $E$ is defined and two definitions,
geometric and topological, of the Euler class of $E$ are given. The proof
of Theorem 1 follows by identifying the geometric and topological Euler
classes (Theorem 4). Theorem 4 is proved in Section
3. We first prove the theorem  for
2-plane   bundles $E$  via an explicit geometric
construction of the Thom class of $E$. This is extended to direct sums of
plane bundles by using elementary properties of the Euler classes. We
 note that Mathai and Quillen [5] have given a geometric
construction of the Thom class for a vector bundle of arbitrary even
rank. The construction in [5] is based upon equivariant differential
forms and the Chern-Weil homomorphism. 

\kern 2 pt

 I would like to thank Harley Flanders,
Steve Rosenberg, Dan Dreibelbis, and David Groisser for 
their help in the preparation of this paper.  I am  also indebted to the
referee for helpful suggestions that resulted in the improvement of the
manuscript.

\goodbreak

\k

\k

\n {\bf 2. Euler Classes}

\kern 2 pt

\n Let $E$ denote a real oriented Riemannian vector 
bundle of even rank $p = 2q$ over a compact connected smooth manifold $M$
(note that $M$ itself does not need to be Riemannian here). In this section, we
associate two De Rham cohomology classes to $E$. 

\k

\n {\bf Definition} A {\it connection} on $E$ is a globally defined
map
$\nabla: \Gamma(E) \mapsto \Gamma(E \otimes T^*(M))$ such that for  
$f \in C^\infty(M)$ and $X \in \Gamma (E),$
$$\nabla(fX) = df\otimes X +f\nabla X. \eqno (2) $$

We say that $\nabla$ is
{\it metric}  if, for all $X$ and $Y \in \Gamma(E)$, we have
$$d<X, Y> = <\nabla X, Y> + <X, \nabla Y>. \eqno (3)$$

\n Suppose that $\nabla$ is a metric connection on $E$. For
every (locally defined) orthonormal frame $e = \{e_1, \dots,
e_p\}^t$ of $E$, we define a $p \times p$ matrix of {\it
connection 1-forms}
$\omega = [\omega_{ij}]$ on
$M$ by the relations
$$\nabla e_i = \sum_j \omega_{ij}e_j \kern 20 pt (\nabla e = \omega e)$$
and a corresponding $p \times p$ matrix of {\it curvature 2-forms}
$\Omega$ by
$$\Omega = d\omega - \omega^2$$
where $d$ denotes exterior derivative and multiplication of matrices of
1-forms is defined in the usual way, using wedge product to multiply
the entries. It is clear that the metric compatability condition (3)
implies that both
$\omega$ and $\Omega$ are skew-symmetric mtrices.
Suppose now that $e$ and $f$ are two positively oriented orthonormal frames of
$E$
defined over intersecting neighborhoods $U$ and $V$ of $M$   with
connection 1-forms and curvature 2-forms
$\omega_e,
\Omega_e$ and
$\omega_f, \Omega_f$ respectively. Then there exists an 
SO(p)-valued  function
$A$ on
$U \cap V$ such that $f = Ae$. Condition (2) implies the
following transformation laws for $\omega$ and $\Omega$:  
$$\omega_f = (dA)A^{-1} + A\omega_eA^{-1} \eqno(4) $$
$$\Omega_f = A\Omega_eA^{-1}. \eqno(5)$$

\k

\n {\it The Pfaffian}: There exists a map $Pf$ 
$: so(p)\mapsto {\bf R}$ (where $so(p)$ denotes the set of real $p\times
p$ skew-symmetric matrices) such that
$Pf(M)$ is a homogeneous polynomial of degree $q = p/2$ in the entries of
$M$, characterized by the following properties:

\kern 10 pt

\n (i) $Pf \pmatrix {0 &\lambda\cr -\lambda &0} = \lambda$.

\n (ii) $Pf \pmatrix {A &0 \cr 0 &B} = 
Pf(A)Pf(B)$, if $A$ and $B$ are skew-symmetric square blocks.

\n (iii) $Pf(AMA^{-1}) = Pf(M), A \in SO(p).$

\k

\n Since the set of even-degree differential form a
commutative ring with
$\wedge$ as multiplication we can define $Pf(\Omega_e)$ where, as
before, $\Omega_e$ is the curvature matrix corresponding to a positively
oriented
orthonormal frame $e$ of $E$ defined over a neighborhood $U$ of $M$. By
(5) and (iii) above, it follows that if $Pf(\Omega_f)$ is another such
expression defined over $V$  with $U \cap V \ne \phi$, then on $U \cap V$   we
have  
$Pf(\Omega_e) = Pf(\Omega_f)$. 
Thus $Pf(\Omega_e)$ extends to a {\it globally}
defined $p$-form on $M$. We denote this by $Pf(\Omega)$.

 \begin{thm}  $Pf(\Omega)$ is closed. Furthermore,
$[Pf(\Omega)]
\in H^p(M)$ is independent of the particular
choice of connection on $E$ used in its construction.
\end{thm}

\n A proof of this result can be found in Rosenberg [R].

\k

\n Define the {\it geometric Euler class} $e_g$ of $E$ by
$$e_g \equiv [(-1/2\pi )^q Pf(\Omega)] \in H^p(M).$$
Now consider $E$ as a (non-compact) manifold of dimension
$n + p$ and denote by
$H_c^p(E)$ the p-th cohomology group of $E$ defined by
compactly supported forms.

\begin{thm}[Thom Isomorphism Theorem] There exists a unique
element
$u
\in H_c^p(E)$, known as the Thom class of $E$, such that for each fiber
$E_x$ of $E$
$$\int_{E_x} u = 1. \eqno(6)$$
\end{thm}

\n See Bott and Tu [2] for a proof.

\k

\n Let $i:M \mapsto E$ denote embedding as the
0-section.  The {\it topological Euler class $e_t$ of $E$} is defined  by
$$e_t \equiv i^*(u) \in H^p(M).$$
Uniqueness of the Thom class implies that $e_t$ is well-defined on the
isomorphism class of
$E$, i.e. it is a  {\it topological invariant}. 
In the case when
$p = n$, i.e. the rank of $E$ and the dimension of $M$ coincide, 
the {\it Euler characteristic} of E is defined by
$$\chi(E) \equiv \int_M e_t.$$

\k

It is an easy exercise to show that  $e_g$
and
$e_t$ (both of which we denote for now by $e(E))$ share the following two
fundamental properties:

\k

\n (i) {\it Whitney duality}: Let $E = E_1\oplus \dots \oplus
E_r$ be a direct sum of oriented even-dimensional bundles. Then 

$$ e(E) = e(E_1)\wedge\dots\wedge(E_r)$$

\n (ii) {\it Naturality}: Let $f: N \mapsto M$ be a $C^\infty$ map and
consider the pull-back bundle $f^*(E)$ over $N$

\kern -10 pt

$$f^*( E) \kern 25 pt E$$
$$\kern 10 pt \big\downarrow \kern 35 pt  \big\downarrow   $$
$$\kern 12 pt N \kern 5 pt \longrightarrow \kern 5 pt M$$

\kern -20  pt

$$\kern 8 pt f$$

\n Then $e(f^*(E)) = f^*(e(E)).$

\k

\n Our main result is

\n \begin{thm}   Let $E$ denote a real orientable Riemannian vector
bundle of even rank $p = 2q$ over a compact orientable manifold $M$ of
dimension $n$. Then
$$e_t = e_g \in H^p(M). \eqno (7)$$
\end{thm}

 Note that if $p = n = 2k$, then integrating each side of (7) over $M$
gives
 Theorem 1.
The Gauss-Bonnet-Chern Theorem is obtained from Theorem 1 by taking 
 $E$ to be  the
tangent bundle of an orientable Riemannian manifold 
$M$, endowed with the Levi-Civita connection.

\k

\k

\n{\bf 3. Proof of Theorem 4}

\kern 2 pt

\n We first prove the theorem for the case where $E$ is a bundle of rank
2, equipped  with a metric connection $\nabla$. The idea is to give an
explicit construction of the Thom class $u$ of $E$ in terms of $\nabla$.

Let 
$\{e_1, e_2\}$  be a (locally defined) positively oriented
orthonormal frame for
$E$ and $\omega_e$ the upper off-diagonal entry of the
corresponding connection 1-form with respect to this frame. Let
$v_1, v_2$ denote the components of $E$-vectors with respect to the frame
$\{e_1, e_2\}$ and
$r$ the radial distance in any fiber $E_x$. Finally,  $\rho$ and
$\gamma$ will denote smooth real-valued functions defined on $[0, \infty)$   of
compact
support and $c$ a  constant, all to be  determined
in the course of the proof. Consider the following locally defined 2-form
on
$E$, which we introduce as a ``template" for constructing the Thom class 
$$u \equiv c\{\rho(r^2)dv_1dv_2 +\rho(r^2)rdr 
\pi^*(\omega_e) + \pi^*(d\omega)\gamma(r)\}, \eqno(8)$$
where $\pi: E \mapsto M$ is the projection map. Note that $r$ and
$d\omega$ are intrinsic objects.  Changing to another orthonormal frame
$\tilde e$ with the same orientation as $e$, related to $e$ by a
conterclockwise rotation of $\theta$, yields 
$$ d\tilde v_1d\tilde v_2 = dv_1dv_2 -
rdr\pi^*(d\theta).$$ However, the transformation law
(4) implies 
$$\tilde \omega_e = \omega_e + d\theta.$$
Substituting into (8), we obtain
$$u = c\{\rho(r^2)d{\tilde v}_1d{\tilde v}_2
+ \rho(r^2)rdr 
\pi^*({\tilde \omega}_e) + \pi^*(d\omega)\gamma(r)\},$$ 
i.e. $u$ is {\it
invariantly defined} (it was this realization that motivated us to
choose (8) as a general form for $u$).  We shall
define $\rho, \gamma$, and $c$ so that $u$ is a representative of the
 Thom class of $E$. First choose
$\gamma$ so that $\gamma(0) = 1$. Since $rdr = v_1dv_1 + v_2dv_2$ and 
$\pi \circ i = id$, applying $i^*$ to (8) yields
$$i^*(u) = cd\omega \eqno (9)$$
Taking the exterior derivative in (8) gives
$$du = c\pi^*(d\omega)dr\{(\gamma'(r) -
\rho(r^2)r\}.$$ Thus a sufficient condition for $u$
to be closed is 
$\gamma'(r) = \rho(r^2)r.$
Together with the condition  
$\gamma(0) = 1$, this implies
$$\gamma(r) = 1 + \int_0^r\rho(s^2)sds.$$
We now choose $\rho$ to be any smooth function with support contained in
$(0, 1)$ such that
 $$\int_0 ^1 \rho(s^2)sds = -1$$ 
and extend $\rho$ and
$\gamma$ by defining them to  be 0 on $[1, \infty)$.
It follows from  (8) that
$$\int_{E_x} u 
= c\int_0^\infty\int_0^\infty \rho(r^2) dv_1dv_2$$ 
$$= c\int_0^{2\pi}d\phi
\int_0^{\infty}\rho(r^2)rdr$$ 
$$= c\int_0^{2\pi}d\phi
\int_0^{1}\rho(r^2)rdr = -2\pi c. $$
Thus chosing $c = -1/2\pi$ ensures that  $u$  
represents the Thom class. With this choice
of $c$, (9) implies 
$$e_t = i^*(u) = e_g $$
and the theorem is proved for the plane bundle case.

Suppose now that $E = E_1\oplus \dots\oplus E_q$ is a sum of
oriented plane 
bundles.
Let $E_1, \dots E_q$ have curvature 2-forms with upper off-diagonal
entries
$d\omega_1,
\dots, d\omega_q$ with respect to (any) metric  connections.
Extend these connections (by direct sum)  to
a connection on $E$. With    the ``block
multiplicative" property of the  Pfaffian we have that
$$pf(\Omega) = d\omega_1\wedge \dots \wedge d\omega_q. \eqno (10)$$
Let $E_1, \dots E_q$ have topological Euler classes
$e_t^1,
\dots, e_t^q $. Using (10), the theorem for 2-plane
bundles, and Whitney duality, we obtain
$$e_g(E) \equiv
 (-2\pi)^{-q}Pf(\Omega)  $$  
$$=  (-d\omega_1/2\pi)\dots (-d\omega_q/2\pi)$$
$$\equiv e_g^1e_g^2\dots e_g^q $$
$$= e_t^1e_t^2\dots e_t^q $$
$$ = e_t(E).$$
Thus the theorem holds in this case.

We complete the proof by appealing to the following result 
(a proof of which can be
found in [7, Page 196]).

\begin{thm} [Splitting Principle] Let $E$ denote a real,
orientable, even dimensional vector bundle over a manifold $M$. Then
there exists a manifold $N$ and a map $g: N\mapsto M$ such that

(i) $g^*: H^*(M)\mapsto H^*(N)$ is a monomorphism (i.e. is
injective). \hfil\break 
(ii) $g^*(E)$ is a sum of orientable plane 
bundles.
\end{thm}

\n Suppose now $E$ is an arbitrary orientable vector bundle over $M$.
Applying the above result with $E$ and using the already established
coincidence of $e_g$ and $e_t$ on sums of orientable plane bundles and
the naturality of $e_g$ and $e_t$, we have
$$
g^*(e_g(E))= e_g(g^*(E))=e_t(g^*(E))=g^*(e_t(E)).
$$
The result now follows from the injectivity property of the map $g^*$.

\k

\k

\centerline {\bf References }

\kern 3 pt

\n [1] J. M. Bismut. The Atiyah-Singer theorems for classical eliptic
operators: a probabilistic approach, I; the index theorem. {\it J. Funct.
Anal.} {\bf 57} (1984), 56-99

\n [2] R. Bott and L. Tu, {\it Differential Forms in Algebraic
Topology}, Graduate Texts in Mathematics vol. 82, Springer-Verlag,
Berlin, New York, 1982.

\n [3] S.-S. Chern, A simple intrinsic proof of the
Gauss-Bonnet formula for closed Riemannian manifolds, {\it Ann. of
Math.} (2) {\bf 45} (1944) 747-752.

\n [4] E. P. Hsu, Stochastic local Gauss-Bonnet-Chern theorem,
{\it  Jour. Theoret. Probab.} {\bf 10}, No. 4,  (1997) 819-834.

\n [5] V. Mathai and D. Quillen, Superconnections, Thom Classes, and
equivariant differential forms, {\it Topology} {\bf 25} (1986) 85-110.

\n [6]  S. Rosenberg, {\it The Laplacian on a Riemannian Manifold},
Cambridge University 
Press,  Cambridge, 1997.

\n [7] P. Shanahan, {\it The Atiyah-Singer Index Theorem: An
Introduction}, Lecture Notes in Mathematics, vol. 638,
Springer-Verlag, Berlin, New York, 1978.

\end{document}